\documentstyle[amssymb,amsmath,amsthm,draft,epsf,12pt]{article}
\scrollmode
\theoremstyle{remark}

\theoremstyle{definition}

 \makeatother

\begin{document}

{\large

\begin{center}

{\LARGE \bf  Hyperbolicity and solvability for linear systems on time scales}
\bigskip

{\bf Sergey Kryzhevich${}^{a,b}$ \footnote{Email address:
kryzhevicz@gmail.com}}

\bigskip

{\noindent\small
${}^a$ Faculty of Mathematics and Mechanics, Saint-Petersburg State University,\\
28, Universitetskiy pr., Peterhof, Saint-Petersburg, 198503, Russia;\\
${}^b$ University of Nova Gorica, Vipavska,13, Nova Gorica, SI-5000, Slovenia}

\bigskip

\end{center}}

\noindent\textbf{Abstract}

We believe that the difference between time scale systems and ordinary differential equations is not as big as people use to think. We consider linear operators that correspond to linear dynamic systems on time scales. We study solvability of these operators in ${\mathbb L}^\infty$. For ordinary differential equations such solvability is equivalent to hyperbolicity of the considered linear system. Using this approach and transformations of the time variable, we spread the concept of hyperbolicity to time scale dynamics. We provide some analogs of well-known facts of Hyperbolic Systems Theory, e.g. the Lyapunov--Perron theorem on stable manifold.

\medskip


\section{Introduction.} 
Time scale systems play an important role in modern dynamics as they stand between discrete and continuous ones. For applications, they could be used for modelling strongly nonlinear phenomena e.g. impacts. There are hundreds of books and papers, devoted to time scale dynamics (see [1-7],  [9], [11-13], [15-20], [23,24], [28,29] and references therein, the list is still incomplete). The main obstacle to study such systems is that they are in principle non-autonomous unless the time scale is periodic. 

Here, we are mostly interested in stability of solutions of time scale systems. There were two principal approaches. One is related to Grobman--Bellman, Bihari and other similar estimates [6], [11], [12,13], [16], [23,24], [28,29], see also [8] for the classical case of ordinary differential equations. Another powerful tool is the so-called second Lyapunov method, related to constructing so-called Lyapunov functions ([2], [4-7], [12], [15-19], [23], [28], see also the classical book [21] for origins). However, in the ODE theory there is the third approach, the so-called first or direct Lyapunov method  [10,21]. Unlike implicit methods, listed above, this method allows to construct  bounded solutions and even invariant manifolds as limits of so-called successive approximations. The main aim of this paper is to generalise this approach, developed for non-autonomous ODEs, to the case of time scale dynamics. We study solvability of operators, corresponding to linear systems, we give analogs of classical result of hyperbolic theory:  existence of bounded solutions for almost linear systems, Lyapunov--Perron theorem on invariant manifolds, etc. 

A similar approach was developed in papers [12] and [28], the principal difference of our approach is that we study equivalences between time scale equations and ODEs. This leads to different results.
The key point is that many linear time scale systems can be represented as reductions of linear systems of ordinary differential equations and solvability of linear time scale operator follows from one of the differential operator. 

In our paper, we always operate with the so-called $\Delta$ -- derivatives, the case of $\nabla$ -- derivatives may be considered similarly. Studying the case of solvability of linear differential operators (and of the time scale ones), we always concentrate on results, related to hyperbolicity (exponential dichotomy) of the corresponding ODE systems. We could also consider the so-called regularity of linear systems or one of its generalisations instead (this would give solvability in the space of exponentially decaying solutions). However, we prefer to postpone this activity for the future.  In this paper, we consider both systems on time scales and ordinary differential equations. We distinguish these two cases by the following formalism: solutions related to time scales are highlighted in bold.  We use standard notions $B(\varepsilon,x)$ for $\varepsilon$ -- ball, centred in $x$ and $|\cdot|$ for the Euclidean norm.

\section{Dynamic systems on time scales}

\noindent\textbf{Definition 2.1.} Let the \emph{time scale} be an unbounded closed subset of $[0,+\infty)$. 

Let $\mathbb T$ be a time scale. Without loss of generality, we always assume that $0\in {\mathbb T}$. 

\noindent\textbf{Definition 2.2.}  Given a $t_0\in {\mathbb T},$ we denote 
$\sigma(t_0):=\inf \{t\in {\mathbb T}: t>t_0\},$ $\mu(t_0):=\sigma(t_0)-t_0$. Such $\mu(t_0)$ is called \emph{graininess} function. We say that $t_0$ is \emph{right-dense} if $\mu(t_0)=0$ and \emph{right-scattered} otherwise. We say that a function $f: {\mathbb T}\to {\mathbb R}$ is \emph{$rd$-continuous} if it is continuous at all right-dense points and left continuous at all left-dense points. 

\noindent\textbf{Definition 2.3.} The function $f : {\mathbb T} \to {\mathbb R}$ is called 
$\Delta$-\emph{differentiable} at a point  $t \in {\mathbb T}$ if there exists 
$\gamma \in {\mathbb R}$ such that for any $\varepsilon > 0$ there exists a
neighborhood ${\mathbf W}\subset {\mathbb T}$ of $t$ satisfying
$$|[f(\sigma(t)) - f(s)] - \gamma[\sigma(t) - s]| \le  \varepsilon|\sigma(t) - s|$$ 
for all $s \in {\mathbf W}$. In this case, we write $f^\Delta(t) = \gamma$. 

When 
${\mathbb T} = {\mathbb R},$ 
$x^\Delta(t) = \dot x(t)$. When  ${\mathbb T} = {\mathbb Z},$ $x^\Delta(n)$ is the standard forward difference operator 
$x(n + 1) - x(n)$.\medskip

\noindent\textbf{Definition 2.4.} If $F^\Delta(t) = f(t),$ $t \in {\mathbb T},$ then $F$ is a 
$\Delta$-\emph{antiderivative} of $f,$ and the Cauchy $\Delta$-\emph{integral} is given by the formula
$$
\int_\tau^s f(t)\Delta t =F(s)-F(\tau) \qquad \mbox{for all} \quad s,\tau\in{\mathbb T}.
$$

Similarly, we may differentiate and integrate vector and matrix-valued functions. 

\noindent\textbf{Definition 2.5.} A function $p : {\mathbb T} \to {\mathbb R}$ is called \emph{regressive} provided  
$1+\mu(t)p(t)\neq 0$ for all $t\in {\mathbb T}$ and \emph{positively regressive} if $1+\mu(t)p(t)> 0$ for all $t\in {\mathbb T}$. 
The set of all regressive and rd-continuous functions is denoted by 
${\mathcal R} = {\mathcal R}({\mathbb T}, {\mathbb R})$. The set of all positively regressive and rd-continuous function is denoted by ${\mathcal R}^+$.

\noindent\textbf{Definition 2.6.} A real non-degenerate matrix $A$ is called \emph{positive} if one of following three equivalent conditions is satisfied:
\begin{enumerate}
\item there is a real matrix $B$ such that $A=\exp(B);$
\item there is a real matrix $C$ such that $A=C^2;$
\item for any negative value $\lambda$ and for any $k\in {\mathbb N}$ the number of entries of the $k\times k$ block 
$$
B_\lambda=\begin{pmatrix} 
\lambda & 0 & 0 &\ldots & 0 \\
1 & \lambda &  0 & \ldots & 0 \\
0 &1 & \lambda & \ldots & 0 \\
\ldots &\ldots & \ldots & \ldots & \ldots \\
0 &\ldots & 0 & 1 & \lambda 
\end{pmatrix}
$$ 
in the Jordan normal form of the matrix $A$ is even (that can be $0,$ of course).
\end{enumerate}
Particularly, the positivity implies (but is not equivalent to) the fact that $\det A>0$. 

Now, we introduce a result from linear algebra. Let $M_{n,n}({\mathbb R})$ (or
$M_{n,n}({\mathbb C})$) be the class of all real (or, respectively, complex) $n\times n$ matrices.

\noindent\textbf{Proposition 2.7.} \emph{There exists a function $\log$ from the set of all non-degenerate $n\times n$ matrices such that the following holds.
\begin{enumerate}
\item $B=\log A$ implies $A=\exp B;$
\item this function is measurable and bounded on any set $\{A: \max (|A|, |A^{-1}|)\le R\},$ $R>0;$
\item if $A$ is positive, then $\log A$ is real-valued.
\end{enumerate}}

The construction of such logarithm is described in [14, Chapter VIII, \S8].

\noindent\textbf{Definition 2.8.} A matrix-valued mapping ${\mathbf A} : {\mathbb T} \to {M}_{n,n}({\mathbb R})$ is called \emph{regressive} if for each 
$t \in {\mathbb T}$ the $n\times n$ matrix $E_n+\mu(t){\mathbf A}(t)$ is invertible, and \emph{uniformly regressive} if in addition the matrix-valued function 
$(E_n+\mu(t){\mathbf A}(t))^{-1}$ is bounded. Here $E_n$ is the unit matrix. We say that the matrix-valued function $\mathbf A$ is positively regressive if all matrices $E_n+\mu(t){\mathbf A}(t)$ are positive. \medskip 

\noindent\textbf{Definition 2.9.} We say that a time scale $\mathbb T$ is \emph{syndetic} if $\sup \{\sigma (t):t\in {\mathbb T}\}<+\infty$ or, in other words, gaps of the time scale are bounded.

We introduce a notion $[t]_{\mathbb T}= \max\{\tau \in {\mathbb T}: \tau\le t\}$. Clearly, $[t]_{\mathbb T}\le t$ and $[t]_{\mathbb T}= t$ if and only if $t\in {\mathbb T}$.

\section{Solvability of linear non-homogenous systems.}

Consider a time scale $\mathbb T$ and an $rd$-continuous matrix-valued function $A:{\mathbb T} \to {\mathbb R}^n$. 

We study a linear system
$$x^\Delta ={\mathbf A}(t) x+{\mathbf f}(t) \eqno (3.1) $$
and the corresponding homogeneous system 
$$x^\Delta ={\mathbf A}(t) x. \eqno (3.2)$$
Here $\mathbf A$ is a bounded uniformly regressive $rd$-continuous matrix-valued function, $\mathbf f$ is a bounded $rd$-continuous vector function. We are interested when systems $(3.1)$ have bounded solutions for all admissible right-hand sides $\mathbf f$. We recall a notion from the theory of linear systems of ordinary differential equations. Given a linear system 
$$\dot x=A(s) x \eqno (3.3)$$
of ordinary differential equations, we consider the Cauchy matrix 
$\Phi_A(t,\tau)=\Phi_A(t)\Phi_A^{-1}(\tau)$. 

\noindent\textbf{Definition 3.1.} A linear system $(3.3)$ is called \emph{hyperbolic} if for any $t\in [0,\infty)$ there exists linear spaces $U^+(t)$ and $U^-(t)$ called stable and unstable spaces respectively and positive values $C$ and $\lambda_0$ such that
\begin{enumerate}
\item $U^+(t)\oplus U^-(t)={\mathbb R}^n;$
\item $\Phi_A(t,\tau)U^\pm(\tau)=U^\pm(t);$
\item $|\Phi_A(t,\tau) x_0| \le C\exp(-\lambda_0(t-\tau)) |x_0|$ for all $t>\tau,$ $x_0\in U^+(\tau);$
\item $|\Phi_A(t,\tau) x_0| \le C\exp(\lambda_0(t-\tau)) |x_0|$ for all $t<\tau,$ $x_0\in U^-(\tau)$.
\end{enumerate}  

Many examples of hyperbolic systems, e.g., for linear systems with constant matrices may be constructed, using approaches of the paper [28].

If a continuous function $f:[0,+\infty)\to {\mathbb R}^n$ is bounded, the system  
$$
\dot x=A(s)x+f(s)
\eqno (3.4)$$
has a bounded solution $\varphi:={\mathcal L} f,$ where
$$\varphi(s)=\int_0^s \Phi_A(s,\tau) \Pi^+(\tau) f(\tau)\, d\tau - \int_t^\infty \Phi_A(t,\tau) \Pi^-(\tau) f(\tau)\, d\tau. \eqno (3.5)$$
Here $\Pi^+(s)$ and $\Pi^-(s)$ are linear projector operators on the stable and the unstable spaces respectively such that $\Pi^+(s)x+\Pi^-(s)x \equiv x$. A similar fact is true for exponentially decaying right hand sides. There exists a $\lambda_1>0$ and $K>0$ such that for any $\lambda\in [0,\lambda_1]$ if 
$|f(t)|\le C\exp(-\lambda t),$ then $|{\mathcal L}f (t)|\le K_\lambda C\exp(-\lambda t)$. Actually, we may take any $\lambda_1 \in (0,\lambda_0)$. The inverse statement is also true (see [22,26] and also [30] for discrete case).

\noindent\textbf{Theorem 3.2.} (Pliss--Maizel Theorem) \emph{If system $(3.4),$ defined on $[0,\+\infty)$ has a bounded solution for any bounded function $f,$ the corresponding system $(3.3)$ is hyperbolic.}

\section{Transformation of the time variable.}

Given a time scale $\mathbb T,$ we define the function $s:{\mathbb R}\to {\mathbb R}$:
$$s(t)=\int_0^t \dfrac{\log (1+\mu ([t]_{\mathbb T}))}{\mu ([t]_{\mathbb T})}\, d\tau.$$  Observe that $s(0)=0$. The following statement is evident. 

\noindent\textbf{Lemma 4.1.} \emph{For any time scale $\mathbb T$ the function $s(t)$ is strictly increasing and unbounded;
$$\limsup_{t\to +\infty} s(t)/t\le 1.$$
If the time scale is syndetic, we also have 
$\liminf\limits_{t\to+\infty} s(t)/t> 0.$} 

Let $\Psi_{\mathbf A}(s,0)$ be a fundamental matrix of the time scale system $(3.2),$ such that $\Psi(0)=E_n$. We prove the following statement.

\noindent\textbf{Lemma 4.2.}  \emph{Given an rd-continuous uniformly regressive matrix 
${\mathbf A}:{\mathbb T}\to {\mathbf M}_{n,n}$ there exists a piece-wise continuous complex matrix-valued function $A: [0,+\infty)\to {\mathbf M}_{n,n}$ such that for 
$$
\Phi_A(s(t),0)=\Psi_{\mathbf A}(t,0)
\eqno (4.1)$$
for all $t\in {\mathbb T}$. If $\sup |{\mathbf A}(t)|<+\infty$ and $\sup |{\mathbf A}^{-1}(t)|<+\infty,$ then $\sup |A(t)|<+\infty$. If $\mathbf A$ is uniformly positively regressive, then the matrix 
$A$ can be taken real.}

\noindent\textbf{Proof.} We set $A(s(t))={\mathbf A}(t)$ for all $t\in {\mathbb T}$. For $t\notin {\mathbb T},$ we set 
$$A(s(t)) = \dfrac{\log [E_n+ \mu ([t]_{\mathbb T}) {\mathbf A}(([t]_{\mathbb T})]}{\log(1+\mu([t]_{\mathbb T}))}.$$
By choice of the function $s(t)$ equality $(4.1)$ is fulfilled. Evidently, 
$$\dfrac{\log(E_n+\mu A)}{\log(1+\mu)}\to A$$ 
as $\mu\to 0$ uniformly on compact sets of matrices $A$. On the other hand,
$$
\lim_{\mu\to+\infty} \dfrac{\log(E_n+\mu A)}{\log(1+\mu)}=
\lim_{\mu\to+\infty} \dfrac{\log \mu E_n+ \log(A+\mu^{-1}E_n)}{\log\mu}=E_n
\eqno (4.2)$$
for any non-degenerate matrix $A$ and the limit is uniform on all compact subsets of $M_{n.n}$ that do not contain degenerate matrices.
$\square$

\noindent\textbf{Definition 4.3.} Consider a time-scale system $(3.2)$ with a uniformly regressive matrix ${\mathbf A}(t)$ such that 
$\sup\{|{\mathbf A}(t)|+|{\mathbf A}^{-1}(t)|:t \in {\mathbb T}\}<+\infty$.
We call it \emph{hyperbolic} if the corresponding system of ordinary differential equations $(3.3)$ is hyperbolic.

For hyperbolic time scale systems, we may take stable and unstable spaces $U^\pm(t)$ (same as for the corresponding systems of ordinary differential equations). 

\noindent\textbf{Proposition 4.4.} \emph{If $(3.2)$ is hyperbolic, there exist constants $C,\lambda>0$ such that
\begin{enumerate}
\item $|\Psi_{\mathbf A}(t,t_0)x_0| \le C |x_0| \exp(-\lambda (s(t)-s(t_0)))$ for all $t,t_0 \in {\mathbb T},$ $t\ge t_0,$ $x_0\in U^{+}(t_0);$
\item $|\Psi_{\mathbf A}(t,t_0)x_0| \le C |x_0| \exp(\lambda(s(t)-s(t_0)))$ for all $t,t_0 \in {\mathbb T},$ $t\le t_0,$ $x_0\in U^{-}(t_0)$.
\end{enumerate}}

Particularly, this statement implies that $\Psi_{\mathbf A}(t,t_0)x_0\to 0$ as $t\to+\infty$ if $x_0\in U^+(t_0)$ and $\Psi_{\mathbf A}(t,t_0)x_0\to \infty$ as $t\to+\infty$ if $x_0\in U^+(t_0) \setminus \{0\}$.

\noindent\textbf{Remark 4.5.} It follows from $(4.2)$ that for any hyperbolic system $(3.2)$ on a time scale ${\mathbb T}$ the following dichotomy takes place: either the time scale is syndetic or the system $(3.2)$ is \emph{unstable hyperbolic} i.e. $U^-(t) \equiv {\mathbb R}^n$.

\section{Transformation of the right hand side.}

Now, we consider a system $(3.2)$ on a time scale ${\mathbb T}$. We fix the corresponding transformation $s(\cdot)$ of the time variable and the corresponding system $(3.3)$ of ordinary differential equations. Suppose that the matrix ${\mathbf A}(t)$ is regressive and invertible for all $t$.
Observe that on the time-scale ${\mathbb T}$ there exists a sigma--algebra, engendered from $\mathbb R,$ so we can consider measurable functions on $\mathbb T$. Given a vector function ${\mathbf f}\in {\mathbb L}^\infty ({\mathbb T}\to {\mathbb R}^n),$ we construct a function $f\in {\mathbb L}^\infty ({\mathbb R}\to {\mathbb R}^n)$ such that 
\begin{enumerate}
\item $f|_{\mathbb T}={\mathbf f};$
\item for any $x_0\in {\mathbb R}^n$ and any $t\in {\mathbb T}$
$$
{\mathbf x}(t,0,x_0)=x(s(t),0,x_0),
\eqno (5.1) $$
where ${\mathbf x}(t,0,x_0)$ ($x(s(t),0,x_0)$) is the solution of systems $(3.2)$ (or, respectively $(3.3)$) with initial conditions $x(0)=x_0;$
\item $f|_{(t,\sigma(t))}=\mbox{\rm const}$ for any $t\in {\mathbb T}$.
\end{enumerate}
By $(4.1)$ (Lemma 4.1) it suffices to check $(4.1)$ for $x_0=0$ only. 
Then $(5.1)$ is equivalent to 
$$
\mu(t_0) {\mathbf f}(t_0)=\int_{s_0}^{s_1} \Phi_A(s_1,\tau) f(\tau)\, d \tau.
\eqno (5.2) $$
Here $s_0=s(t_0),$ $s_1=s(\sigma(t_0))$. In our assumptions, setting $f_0:=f|_{(s_0,s_1)},$ we reformulate Eq.\, $(5.2)$ as follows: 
$$
{\mathbf f}(t_0)= \dfrac{A^{-1} (\exp(A(s_0)(s_1-s_0)) -E_n) f_0}{\mu(t_0)}= \log(1+\mu(t_0))
(\log [E_n+\mu(t_0) {\mathbf A}(t_0)])^{-1} {\mathbf A}(t_0) f_0
$$ 
if $\mu(t_0)>0$ or 
$$
f_0=\dfrac{\log[E_n+\mu(t_0){\mathbf A}(t_0)]}{\log(1+\mu(t_0))} {\mathbf A}^{-1}(t_0){\mathbf f}(t_0)
$$
These formulae imply the following statement.

\noindent\textbf{Theorem 5.1.}  \emph{Let the matrix $\mathbf A$ be uniformly regressive with respect to the time scale ${\mathbb T},$ hyperbolic and uniformly bounded together with the inverse matrix ${\mathbf A}^{-1}$. Then, for any function ${\mathbf f}\in {\mathbb L}^\infty ({\mathbb T}\to {\mathbb R}^n)$ the corresponding system $(3.1)$ has a bounded solution.}

In this case, there exists a continuous linear operator 
${\mathbf L}: {\mathbb L}^\infty ({\mathbb T}\to {\mathbb R}^n) \to {\mathbb L}^\infty ({\mathbb T}\to {\mathbb R}^n)$
such that for any ${\mathbf f}\in {\mathbb L}^\infty ({\mathbb T}\to {\mathbb R}^n)$ the function ${\mathbf L} {\mathbf f}$ is a bounded solution of system $(3.1)$. This operator $\mathbf L$ corresponds to the operator $\mathcal L$ that gives a bounded solution for Eq.\, (3.4) and is defined by formula (3.5). Let ${\mathbf K}=\|\mathbf L\|$.

An analog of Pliss--Maizel Theorem is also true for time scale systems.

\noindent\textbf{Theorem 5.2.}  \emph{Let the matrix $\mathbf A$ be uniformly regressive with respect to the time scale ${\mathbb T}$. Suppose that for any ${\mathbf f}\in {\mathbb L}^\infty$ the corresponding system $(3.1)$ has a bounded solution and the time scale is syndetic. Then system $(3.2)$ is hyperbolic on $\mathbb T$.}

\noindent\textbf{Proof.} Suppose that system $(3.3),$ constructed by system $(3.2)$ as demonstrated is not hyperbolic. Then, there exists a bounded right hand side $f$ such that the corresponding system $(3.4)$ does not have any solutions, bounded on $[0,+\infty)$. Since system $(3.4)$ is linear, all coefficients are bounded and the time scale is syndetic, all solutions of $(3.4)$ are unbounded on $\mathbb T$. Consider the function ${\mathbf f}:{\mathbb T}\to {\mathbb R}^n$ such that $f(t)={\mathbf f}(t)$ for all right-dense points $t$ and Eq. (5.2) is satisfied for all right-scattered points. Then all solutions of Eq. (3.1) are unbounded. $\square$

Similarly to what is done for ordinary differential equations, we can give estimates of the operator 
$\mathbf L$ in spaces of "exponentially small"\ functions.

\noindent\textbf{Proposition 5.3.} \emph{Let the matrix $\mathbf A$ be hyperbolic on the time scale $\mathbb T$. Then there exist $K>0$ and $\lambda_1>0$ such that for any $\lambda\in [0,\lambda_0]$ the inequality 
$|{\mathbf f}(t)|\le C\exp(-\lambda s(t))$ $\forall t \in {\mathbb T}$
implies $|{\mathbf Lf} (t)|\le CK \exp(-\lambda s(t)) \quad \forall t \in {\mathbb T}$.}

\section{Conditional stability by first approximation}

We can use the statement of Theorem 4.1 to prove some time scale analogs of famous statements from the theory of hyperbolic ODEs. 

\noindent\textbf{Theorem 6.1.}  \emph{Let the matrix $\mathbf A$ satisfy conditions of Theorem 5.1. Let $r_0>0$ and the continuous function ${\mathbf g}: {\mathbb T}\times B(0,r_0)$ be such that 
\begin{enumerate}
\item $|{\mathbf g}(t,0)|\le \varepsilon$ for any $t\in {\mathbb T};$
\item $|{\mathbf g}(t,x_1)-{\mathbf g}(t,x_2)|\le l |x_1-x_2|$ for any $t\in {\mathbb T},$ $x_{1,2} \in B(0,r_0)$.
\end{enumerate} 
Then given $r_0$ there exist $\varepsilon_0, l_0>0$ such that if $l<l_0,$ $\varepsilon<\varepsilon_0$ there exists a bounded solution ${\mathbf X}(t)$ of the system
$$
x^\Delta={\mathbf A}(t)x+{\mathbf g}(t,x)
\eqno (6.1)$$
such that 
$$
|{\mathbf X}(t)|\le \dfrac{K \varepsilon}{1-Kl}. \eqno (6.2)
$$}

Let $\lambda_0$ be a constant of hyperbolicity of the matrix $\mathbf A$.

\noindent\textbf{Theorem 6.2 (Lyapunov--Perron Theorem).}  \emph{Let the matrix $\mathbf A$ satisfy conditions of Theorem 5.1. Let $r_0>0$ and the continuous function ${\mathbf g}: {\mathbb T}\times B(0,r_0)$ be such that 
\begin{enumerate}
\item ${\mathbf g}(t,0)=0$ for any $t\in {\mathbb T};$
\item $|{\mathbf g}(t,x_1)-{\mathbf g}(t,x_2)|\le l |x_1-x_2|$ for any $t\in {\mathbb T},$ $x_{1,2} \in B(0,r_0)$.
\end{enumerate} 
Then given $r_0,$ $\lambda\in (0,\lambda_0),$ $t_0\in {\mathbb T}$ there exist $D>0$ $l_0>0$ such that if $l<l_0,$  there exists a map 
$h: B(0,r_0) \bigcap U^+(t_0)\to U^-(t_0)$ such that 
\begin{enumerate}
\item $h(0)=0;$
\item $|h(x)-h(y)|\le Dl |x-y|$.
\item If $x_0$ is such that $x_0=y_0+h(y_0)$ for some $y_0,$ then ${\mathbf x}(t,t_0, x_0)$ tends to zero as $t$ goes to infinity (in fact, it tends to zero).
Here ${\mathbf x}(t,t_0, x_0)$ is the solution of system $(5.3)$  with initial conditions ${\mathbf x}(t_0)=x_0$.
\end{enumerate}}

This allows to construct the so-called local stable manifold as the image of the constructed map $h$. By remark 4.6, this result is non-trivial only if the time scale is syndetic. Proofs of Theorems 6.1 and 6.2 are very close to ones of their classical analogs [8,10,21,25-27].

\noindent\textbf{Proof of Theorem 6.1.} 
Consider the equation
$${\mathbf x}(t)={\mathbf L} [{\mathbf g} (\cdot, {\mathbf x})](t). \eqno (6.3)$$

Any solution ${\mathbf X}(t)$ of $(6.3)$ is a bounded solution of the equation 
$x^\Delta={\mathbf A}(t)x+{\mathbf g}(t,X(t))$
and, hence, one of Eq.\, (6.1). Given $r_0,$ we take $\varepsilon_0$ and $l_0$ so small that 
$$Kl_0 \le \dfrac12,\qquad \dfrac{K\varepsilon_0}{1-Kl_0}\le \dfrac{r_0}2.$$

We set ${\mathbf x}^0(t)=0$ for all $t$ and define 
$${\mathbf x}^m(t)={\mathbf L} [{\mathbf g} (\cdot, {\mathbf x}^{m-1})](t). \eqno (6.4)$$
for all $m\in {\mathbb N}$.

\noindent\textbf{Lemma 6.3.} \emph{All approximations ${\mathbf x}^m(t)$ 
$(m\in {\mathbb N} \bigcup \{0\})$ are 
\begin{enumerate}
\item well-defined on $\mathbb T;$ 
\item such that $|{\mathbf x}^m(t)|\le r_0$ for all $t\in {\mathbb T},$ $m\in {\mathbb N};$
\item such that
$$|{\mathbf x}^{m+1}(t) - {\mathbf x}^m(t)| \le {K\varepsilon}(Kl)^m, \qquad t\in {\mathbb T}.
\eqno (6.5)$$
\end{enumerate}}

Observe that,
$$\|{\mathbf x}^1\|=\|{\mathbf x}^1-{\mathbf x}^0\| \le \|{\mathbf L} [{\mathbf g} (\cdot, 0)]\|\le K\varepsilon
\eqno (6.6)
$$
(all norms are considered in ${\mathbb L}^\infty({\mathbb T})$). So, the statement of the lemma is true for $m=0$. 

Proceed by induction from the step $m-1$ to $m$. If $\|{\mathbf x}^{m-1}\|\le r_0,$ the right hand side of Eq. (6.4) is well-defined and the solution ${\mathbf x}^m(t)$ can be found. Inequalities (6.5) considered for all previous steps and (6.6) imply that 
$$\|{\mathbf x}^{m}\|\le \dfrac{K\varepsilon(1-(Kl)^{m+1})}{1-Kl}\le \dfrac{r_0}2. \eqno (6.7)$$

Hence the iteration ${\mathbf x}^{m+1}$ is also well-defined and
$$\|{\mathbf x}^{m+1}-{\mathbf x}^{m}\|=
\|{\mathbf L}[{\mathbf g} (\cdot, {\mathbf x}^m)-{\mathbf g} (\cdot, {\mathbf x}^{m-1})]\|\le Kl \|{\mathbf x}^{m} - {\mathbf x}^{m-1}\|$$
that implies (6.5). $\square$.

So, the iterations ${\mathbf x}^k$ converge uniformly and we may set ${\mathbf X}=\lim {\mathbf x}^m$. Since the function ${\mathbf g}$ is uniformly continuous w.r.t. $x,$ we can proceed to limit in (6.4). So, $\mathbb X$ is a solution of (6.3). Proceeding to limit in Eq. (6.7), we get (6.2) that finishes the proof. $\square$

\noindent\textbf{Proof of Theorem 6.2.} Without loss of generality, we suppose that $t_0=0$. Fix 
$y_0 \in U^+(0)$. Take ${\mathbf x}^0(t,y_0)=0,$ 
${\mathbf x}^1(t,y_0)=\Psi_{\mathbf A}(t,0)y_0$ and set
${\mathbf x}^{m+1}(t,y_0)=\Psi_{\mathbf A}(t,0)y_0+ 
{\mathbf L}[{\mathbf g}(\cdot,{\mathbf x}^{m})](t)$
for all $m\in {\mathbb N}$. By definition, we have $|{\mathbf x}^1(t,y_0)|\le a |y_0| \exp(-\lambda s(t))$. We consider $y_0$ so small that $2a|y_0|<r_0$. We prove the following lemma, similar to Lemma 6.3. 

\noindent\textbf{Lemma 6.4.} \emph{All approximations ${\mathbf x}^m(t)$ 
$(m\in {\mathbb N} \bigcup \{0\})$ are 
\begin{enumerate}
\item well-defined on $\mathbb T;$ 
\item such that 
$$|{\mathbf x}^m(t,y_0)|\le 2a |y_0|\exp(-\lambda s(t))\le r_0, \qquad t\in {\mathbb T},\quad 
m\in {\mathbb N}; \qquad \mbox{and} \eqno (6.8)$$ 
$$|{\mathbf x}^{m+1}(t,y_0) - {\mathbf x}^m(t,y_0)| \le a(Kl)^m |y_0|\exp(-\lambda s(t)), \qquad t\in {\mathbb T}.
\eqno (6.9)$$
\end{enumerate}}

\noindent\textbf{Proof.} Inequalities (6.8) are evident for $m=0$ and $m=1,$ inequality (6.9) is evident for $m=0$. Now we are going to prove the lemma by induction. 

If $|{\mathbf x}^m(t,y_0)|\le r_0$ (see Eq. ), then the next approximation $|{\mathbf x}^{m+1}(t,y_0)|$ is correctly defined. Moreover, 
$|{\mathbf x}^{m+1}(t,y_0) - {\mathbf x}^m(t,y_0)|=|{\mathbf L}[{\mathbf g} (\cdot, {\mathbf x}^m)](t,y_0)-{\mathbf L}[{\mathbf g} (\cdot, {\mathbf x}^{m-1})](t,y_0)|\le$
$$\le Kl |y_0| a(Kl)^m |y_0| \exp(-\lambda s(t))$$
which proves (6.9) for the given $m$. Taking sum of inequalities (6.9) for all previous values of $m$ and taking into account the estimate for ${\mathbf x}^1,$ we get
$$|{\mathbf x}^{m+1}(t,y_0)|\le \left(a+a\dfrac{Kl}{1-Kl}\right)  |y_0|\exp(-\lambda s(t)).$$
If $Kl<1/2,$ this implies (6.8) on the step $m+1$.
$\square$

Now we prove that all iterations ${\mathbf x}^{m}$ are Lipschitz continuous. We set 
$${\mathbf x}^m(t,y_0)={\mathbf x}^1(t,y_0)+{\mathbf z}^m(t,y_0)=\Psi_{\mathbf A}(t,0)y_0+ {\mathbf z}^m(t,y_0).$$

\noindent\textbf{Lemma 6.5.} \emph{All iterations ${\mathbf x}^m(t,y_0)$ and ${\mathbf z}^m(t,y_0)$
are Lipschitz continuous: for any $t\in {\mathbb T},$ $y_{0},y_1$ such that $|y_{0,1}|\le r_0$
$$\begin{array}{c}
|{\mathbf x}^m(t,y_0) - {\mathbf x}^m(t,y_1)|\le 2a \exp(-\lambda s(t))|y_0-y_1|;\\
|{\mathbf z}^m(t,y_0) - {\mathbf z}^m(t,y_1)|\le 2Kal \exp(-\lambda s(t))|y_0-y_1|.
\end{array}\eqno (6.10)$$}

\noindent\textbf{Proof.} For $m=1,$ (6.10) is evident: 
$$|{\mathbf x}^1(t,y_0) - {\mathbf x}^1(t,y_1)|\le a \exp(-\lambda s(t))|y_0-y_1|,$$
${\mathbf z}^1(t,y_0)\equiv 0$. Then, we continue the proof by induction.
 
Let (6.10) be satisfied for a fixed value $m$. We write
$$\begin{array}{c}
|{\mathbf z}^{m+1}(t,y_0) - {\mathbf z}^{m+1}(t,y_1)|=\\
|{\mathbf L}[{\mathbf g}(\cdot,{\mathbf x}^{m+1}(\cdot,y_0)](t)-{\mathbf L}[{\mathbf g}(\cdot,{\mathbf x}^{m+1}(\cdot,y_0)](t)|\le 2aKl \exp(-\lambda s(t))|y_0-y_1|;\\[7pt]
|{\mathbf x}^{m+1}(t,y_0) - {\mathbf x}^{m+1}(t,y_1)|\le |{\mathbf x}^{1}(t,y_0) - {\mathbf x}^{1}(t,y_1)| + |{\mathbf z}^{m+1}(t,y_0) - {\mathbf z}^{m+1}(t,y_1)|\le\\
(a+2Kal)\exp(-\lambda s(t))|y_0-y_1|\le 2a \exp(-\lambda s(t))|y_0-y_1|. \qquad \square
\end{array}$$

By Lemma 6.4, approximations ${\mathbf x}^k(t,y_0)$ converge to 
${\mathbf x}^*(t,y_0)={\mathbf x}^1(t,y_0)+{\mathbf z}^*(t,y_0)$
that is a solution of the equation 
${\mathbf x}(t)=\Psi_{\mathbf A}(t,0)y_0+{\mathbf L}[{\mathbf g}(\cdot,{\mathbf x}(\cdot))](t)$
with initial conditions 
${\mathbf x}(0)=y_0+{\mathbf z}^*(0,y_0)=:y_0+h(y_0)$.
Proceeding to limit in (6.8), we get
$$|{\mathbf x}^*(t,y_0)|\le 2a |y_0|\exp(-\lambda s(t)),$$ the second line of (6.10) implies
$$|h(y_0)-h(y_1)|\le 2aKl |y_0-y_1|. \qquad \square$$ 

Many other analogs of classical results of hyperbolic systems of o.d.e.s may be proved for time scale systems. For example, following the lines of [27, Chapter 1], we can prove that all solutions that start in a small neighbourhood of zero out of the stable manifold, leave this small neighbourhood as time increases. Also, we can prove that for any $r\in {\mathbb N}$ the stable manifold is $C^r$ - smooth provided the function $\mathbf g$ is $C^r$ - smooth w.r.t. $x$.

\medskip 

\noindent\textbf{Acknowledgements.}  The author was partially supported by RFBR grant 15-01-03797-a.

\end{document}